\documentclass{proc-l}

\usepackage{amssymb}

\usepackage{graphicx}

\usepackage{biblatex}

\usepackage{tikz-cd}

\newtheorem{theorem}{Theorem}

\newtheorem*{theorem*}{Theorem}
\newtheorem*{conjecture*}{Conjecture}
\newtheorem*{proposition*}{Proposition 3}

\theoremstyle{definition}

\newtheorem*{example*}{Example}

\theoremstyle{remark}

\usepackage[justification=centering]{caption}

\usepackage{mathrsfs}

\usepackage[scr=rsfs,cal=boondox]{mathalfa}

\bibliography{bibliography.bib}

\begin{document}

\title{Semi-isotopic knots}


\author{Fredric D. Ancel}
\address{Department of Mathematical Sciences, University of Wisconsin-Milwaukee}
\curraddr{}
\email{ancel@uwm.edu}
\thanks{}

\subjclass[2010]{54B17, 57K10, 57K30, 57M30, 57N37}


\dedicatory{}


\begin{abstract}A \emph{knot} is a possibly wild simple closed curve in $S^3$.  A knot $J$ is \emph{semi-isotopic} to a knot $K$ if there is an annulus $A$ in $S^3\times[0,1]$ such that $A\cap(S^3\times\{0,1\})=\partial A=(J\times\{0\})\cup(K\times\{1\})$ and there is a homeomorphism $e:S^1\times[0,1)\rightarrow A-(K\times\{1\})$ such that $e(S^1\times\{t\})\subset S^3\times\{t\}$ for every $t\in[0,1)$.
\begin{theorem*}{Every knot is semi-isotopic to an unknot.}\end{theorem*}  
\end{abstract}

\maketitle

\section{Introduction}
We fix some notation and terminology.  Let $I=[0,1]$.  A \emph{knot} is the image of an embedding $S^1\rightarrow S^3$.  If the composition of this embedding with some homeomorphism of $S^3$ is a piecewise linear embedding, then the knot is tame.  Otherwise, it is wild.  An annulus is a space that is homeomorphic to $S^1\times I$.

Knots J and K are \emph{ambiently isotopic} if there is a level-preserving homeomorphism 
$h:S^3\times I\rightarrow S^3\times I$ such that $h(x,0)=(x,0)$ and $h(J\times\{1\})=K\times\{1\}$.  Note that for such an $h$, $h(J\times\{0\})=J\times\{0\}$ and $h(J\times\{1\})=K\times\{1\}$.  Of course, classical knot theory is the study of ambient isotopy classes of tame knots in $S^3$.

Knots $J$ and $K$ are \emph{(non-ambiently) isotopic} if there is a level-preserving embedding 
$e:J\times I\rightarrow S^3\times I$ such that $e(J\times\{0\})=J\times\{0\}$ and $e(J\times\{1\})=K\times\{1\}$.

Observe that every knot that pierces a tame disk is isotopic to an unknot.  The sequence of pictures in Figure 1 suggests a proof of this observation.

\begin{figure}
    \centering
    \includegraphics[scale=.5]{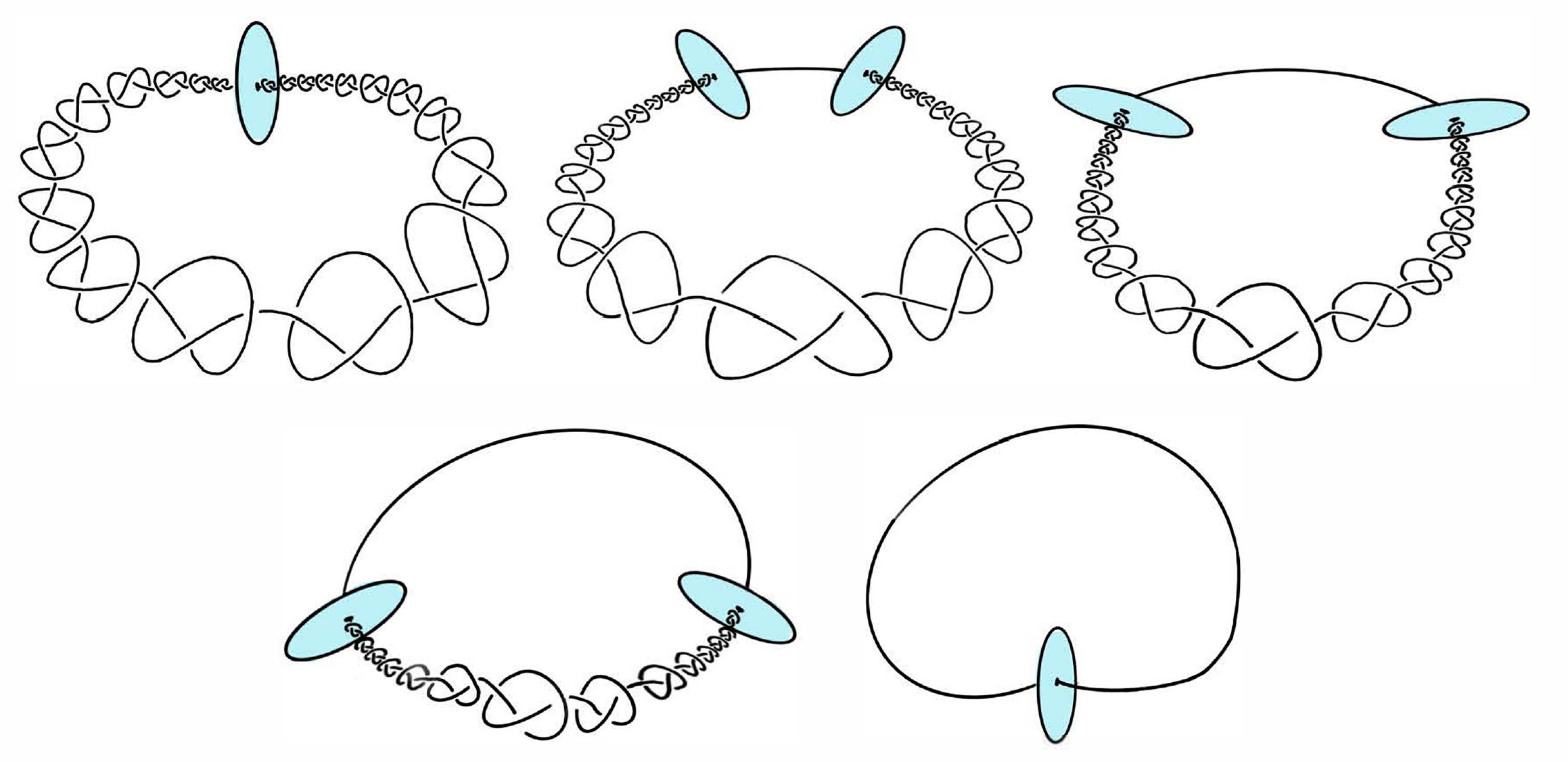}
    \caption{An isotopy of a knot to an unknot}
\end{figure}

The \emph{Bing sling} (Figure 2) \cite{bing1956simple} is a wild knot that pierces no disk.  It is not known whether the Bing sling is isotopic to an unknot.

\begin{figure}
    \centering
    \vspace{1cm}
    \includegraphics[scale=1.9]{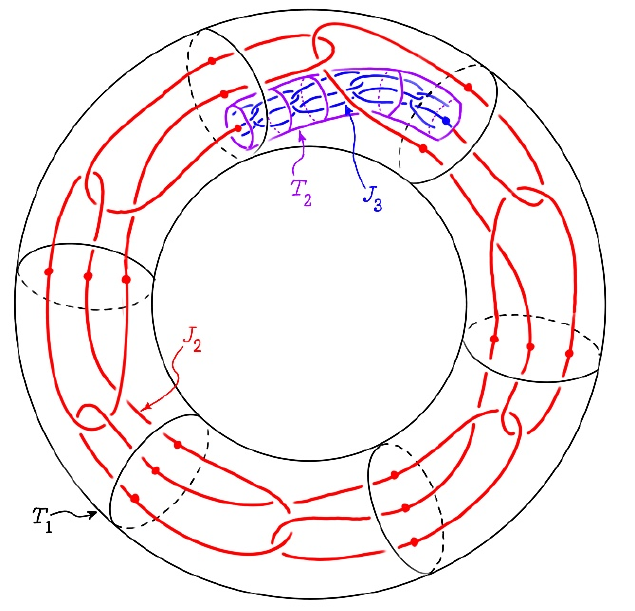}
    \caption{The Bing sling}
\end{figure}

The following conjecture is well-known.
\begin{conjecture*}{Every knot is isotopic to an unknot.}\end{conjecture*}

A knot $J$ is \emph{semi-isotopic} to a knot $K$ if there is an annulus $A$ in $S^3\times I$ such that $\partial A=(J\times\{0\})\cup(K\times\{1\})$ and there is a homeomorphism $e:S^1\times[0,1)\rightarrow A-(K\times\{1\})$ such that $e(S^1\times\{t\})\subset S^3\times\{t\}$ for every $t\in[0,1)$.  Note that $e$ may not extend continuously to homeomorphism from $S^1\times[0,1]$ onto $A$.

The main result of this paper is:
\begin{theorem*}{Every knot is semi-isotopic to an unknot.}\end{theorem*} 

Thus, the Bing sling is semi-isotopic to an unknot.

Knots $J$ and $K$ are \emph{(topologically) concordant} or \emph{I-equivalent} if there is an annulus $A$ in $S^3\times I$ such that $A\cap(S^3\times \{0,1\})=\partial A=(J\times\{0\})\cup(K\times\{1\})$.

Note that: \emph{Isotopic} $\Rightarrow$ \emph{semi-isotopic} $\Rightarrow$ \emph{concordant}.

\begin{example*}{There is a (wild) two-component link in $S^3$ that is not concordant to any PL link!} \cite{melikhov2020topological}\end{example*}

Observe that Melikhov’s striking example shows that the method of proof of the Theorem won’t extend to two-component links, and that the invariants which Melikhov exploits won’t work for single component knots.

The Theorem is proved by applying the two main results of \cite{ancel1996mapping}.  These results are a reinterpretation and formalization of a technique introduced by the topologist C. H. Giffen in the 1960s called \emph{shift-spinning}.  A consequence of Giffen’s method is that the Mazur 4-manifold \cite{mazur1961note} has a disk pseudo-spine.  (A \emph{pseudo-spine} of a compact manifold $M$ is a compact subset $X$ of $int(M)$ such that $M-X$ is homeomorphic to $\partial M\times[0,1)$.)  After R. D. Edwards’ groundbreaking proof that the double suspension of the boundary of the Mazur 4-manifold is homeomorphic to the 5-sphere in 1975, it was observed that had the existence of a disk pseudo-spine in the Mazur 4-manifold been observed before Edwards’ proof, then a proof of Edwards’ theorem using previously known results could have been given.  However, no one made the connection between Giffen’s technique and the existence of a disk pseudo-spine until after Edwards announced his proof.  Shift-spinning, the existence of a disk pseudo-spine in the Mazur 4-manifold and why this leads to the conclusion that the double suspension of the boundary of the Mazur 4-manifold is homeomorphic to the 5-sphere are clearly explained on pages 102-107 of \cite{daverman1986decompositions}.  In \cite{ancel1996mapping} the main results are applied to show that 4-manifolds like the Mazur manifold (constructed by attaching 2-handles to $B^3\times S^1$ have identifiable pseudo-spines.  In particular, all the 4-manifolds constructed by attaching a single 2-handle to $B^3\times S^1$ along a degree one curve have disk pseudo-spines, and therefore have the property that the double suspensions of their boundaries are homeomorphic to $S^5$.  

Other topologists have explored Giffen shift-spinning in various contexts. (See pages 404-409 of \cite{matsumoto1979wild} and pages 15-16 of \cite{hillman2012algebraic}.)

We gratefully acknowledge the artistic hand of Shayna Meyers in the creation of the figures for this article.  We also thank Chris Hruska for a useful conversation.

\section{Mapping swirls}
The central concept of \cite{ancel1996mapping} is the \emph{mapping swirl}.  We recapitulate its definition and state the two main results of \cite{ancel1996mapping} in forms that are convenient for proving the Theorem. 

Let $X$ be a compact metric space.  Identify the \emph{cone} on $X$, $CX$, and the \emph{suspension} of $X$, $\Sigma X$, as the quotient spaces $CX=([0,\infty]\times X)/(\{\infty\}\times X)$ and $\Sigma X=([-\infty,\infty]\times X)/\{\{-\infty\}\}\times X,\{\infty \}\times X \}$.  Let $(t,x)\mapsto tx$ denote either of the quotient maps $[0,\infty]\times X\rightarrow CX$ or $[-\infty,\infty]\times X\rightarrow \Sigma X$.  For each $t\in [0,\infty]$ or $[-\infty,\infty]$, let $tX$ denote the image of the set $\{t\}\times X$ under the appropriate quotient map.  Thus, $\infty X$ denotes the \emph{cone point} of $CX$, and $(-\infty)X$ and $\infty X$ denote the \emph{suspension points} of $\Sigma X$. 

Let $f:X\rightarrow Y$ be a map between compact metric spaces.  Observe that, by exploiting the homeomorphism $x\mapsto(x,f(x))$ from $X$ to the graph of $f$, the \emph{mapping cylinder} of $f$, $Cyl(f)$, can be identified with the subset \[\{(tx,f(x))\in CX\times Y:(t,x)\in[0,\infty)\times X\}\cup({\infty X}\times Y)\] of $CX\times Y$.  
Similarly, the \emph{double mapping cylinder} of $f$, $DblCyl(f)$, which is obtained by 
identifying two copies of $Cyl(f)$ along their bases, can be identified with the subset  
\[\{(tx,f(x))\in\Sigma X\times Y:(t,x)\in(-\infty,\infty)\times X\}\cup(\{(-\infty)X,\infty X\}\times Y)\]
of $\Sigma X\times Y$.  For each $t\in(-\infty,\infty)$, the \emph{t-level} of $DblCyl(f)$ is the set 
\[L(f,t)=(tX\times Y)\cap DblCyl(f)=\{(tx,f(x)):x\in X\}.\]
Furthermore, if $t\in[0,\infty)$, then $L(f,t)$ is also called the \emph{t-level} of $Cyl(f)$.  Note that for each $t\in(-\infty,\infty)$, $x\mapsto(tx,f(x)):X\rightarrow L(f,t)$ is a homeomorphism.  The $\infty$\emph{-level} of $Cyl(f)$ is the set $L(f,\infty)=\{\infty X\}\times Y$.  The $(-\infty)$\emph{-level} and the $\infty$\emph{-level} of $DblCyl(f)$ are the sets $L(f,-\infty)=\{(-\infty)X\}\times Y$ and $L(f,\infty)=\{(\infty)X\}\times Y$.

Let $X$ be a compact metric space and let $f:X\rightarrow S^1$ be a map.  The \emph{mapping swirl} of $f$ is 
the subset \[Swl(f)=\{(tx,e^{2\pi it}f(x))\in CX\times S^1:(t,x)\in[0,\infty)\times X\}\cup(\{\infty X\}\times S^1)\] of $CX\times S^1$.  The \emph{double mapping swirl} of $f$ is the subset $DblSwl(f)=$ \[\{(tx,e^{2\pi it}f(x))\in \Sigma X\times S^1:(t,x)\in(-\infty,\infty)\times X\}\cup(\{(-\infty)X,\infty X\}\times S^1)\] of $\Sigma X\times S^1$.  For each $t\in(-\infty,\infty)$, the \emph{t-level} of $DblSwl(f)$ is the set \[\mathcal{L}(f,t)=(tX\times S^1)\cap DblSwl(f)=\{(tx,e^{2\pi it}f(x)):x\in X\}.\]  Furthermore, if $t\in[0,\infty)$, then $\mathcal{L}(f,t)$ is also called the \emph{t-level} of $Swl(f)$.  Note that for each $t\in(-\infty,\infty)$, $x\mapsto(tx,e^{2\pi it}f(x)):X\rightarrow\mathcal{L}(f,t)$ is a homeomorphism.  The $\infty$\emph{-level} of $Swl(f)$ is the set $\mathcal{L}(f,\infty)=\{\infty X\}\times S^1$.  The $(-\infty)$\emph{-level} and $\infty$\emph{-level} of $DblSwl(f)$ are the sets $\mathcal{L}(f,-\infty)=\{(-\infty)X\}\times S^1$ and $\mathcal{L}(f,\infty)=\{\infty X\}\times S^1$.  For each $x\in X$, the \emph{x-fiber} of $Swl(f)$ is the set
\[\mathcal{F}(f,x)=\{(tx,e^{2\pi it}f(x)):t\in[0,\infty)\},\] and the \emph{x-fiber} of $DblSwl(f)$ is the set
\[\mathcal{DblF}(f,x)=\{(tx,e^{2\pi it}f(x)):t\in(-\infty,\infty)\}.\]

We now state the two main results of \cite{ancel1996mapping}.

\begin{theorem}{If $X$ is a compact metric space and $f,g:X\rightarrow S^1$ are homotopic maps, then 
there is a homeomorphism $\Omega:Swl(f)\rightarrow Swl(g)$ with the following properties.\\
1) $\Omega$ is \emph{fiber-preserving}: for every $x\in X$, $\Omega(\mathcal{F}(f,x))=\mathcal{F}(g,x)$.\\
2) $\Omega$ fixes the $\infty$-level: $\Omega |\mathcal{L}(f,\infty)=id.$}\end{theorem} 

Compare Theorem 1 to the fact that, in general, there is no homeomorphism between mapping cylinders of homotopic maps.

Note that Theorem 1 holds with $S^1$ replaced by any space homeomorphic to $S^1$.  Also note that the conclusions of the theorem imply that for any subset $Y$ of $X$, $\Omega(Swl(f|Y))=Swl(g|Y)$, $\Omega(\mathcal{L}(f|Y,0)=\mathcal{L}(g|Y,0)$ and $\Omega(\mathcal{L}(f,\infty))=\mathcal{L}(g,\infty)$.

\begin{theorem}{If $Y$ is a compact metric space, $n$ is a non-zero integer and $f:X\times S^1\rightarrow S^1$ is the map satisfying $f(y,z)=z^n$, then there is a homeomorphism $\Theta:Cyl(f)\rightarrow Swl(f)$ with the following properties.\\
1) $\Theta$ is level-preserving: for every $t\in[0,\infty]$, $\Theta(L(f,t))=\mathcal{L}(f,t)$.\\
2) $\Omega$ fixes the 0- and $\infty$-levels: $\Theta |L(f,0)\cup L(f,\infty)=id.$}\end{theorem} 

Note that Theorem 2 holds for any map $f:X\rightarrow J$ for which there is a commutative diagram: 
\[
    \begin{tikzcd}
        Y\times S^1   \arrow{rr}{(y,z)\mapsto z^n} \arrow{d}[left]{homeomorphism} && S^1 \arrow{d}{homeomorphism} \\ X \arrow{rr}[below]{f} && J
    \end{tikzcd}
\]

We will need one other elementary fact:

\begin{proposition*}{If $X$ is a compact metric space, $f:X\rightarrow S^1$ a map and 
$\lambda:[0,1)\rightarrow[0,\infty)$ is a homeomorphism, then there is a homeomorphism $\Lambda:X\times [0,1)\rightarrow Swl(f)-\mathcal{L}(f,\infty)$ such that for every $t\in[0,1)$, $\Lambda(X\times\{t\}))=\mathcal{L}(f,\lambda(t))$.}\end{proposition*}
 
Proofs of Theorems 1 and 2 are found in \cite{ancel1996mapping}.  Because we have modified the statements of these theorems for their use in this paper, we will provide outlines of these proofs as well as a proof of Proposition 3 in section 4 below.

\section{The proof of the Theorem}
Let $J$ be a knot.  We will prove that $J$ is semi-isotopic to an unknot.

\textbf{Step 1:}  There is an unknotted solid torus $T$ in $S^3$ such that $J\subset int(T)$ and the inclusion $J\hookrightarrow T$ is a homotopy equivalence.

We can assume $J\subset \mathbb{R}^3=S^3-\{\infty\}$.  Let $p$ and $q$ be distinct points of $J$.  Let $\mathcal{V}$ be an uncountable family of parallel planes in $\mathbb{R}^3$ that separate $p$ from $q$.  Since $J$ has a countable dense subset, $J$ does not contain an uncountable pairwise disjoint collection of non-empty open sets.  Hence, there is a $V\in\mathcal{V}$ such that $J\cap V$ contains no non-empty open subset of $J$.  It follows that $J\cap V$ is a totally disconnected subset of $V$.  Let $J_1$ and $J_2$ be arcs such that $J_1\cup J_2=J$ and $J_1\cap J_2=\partial J_1=\partial J_2=\{p,q\}$.  Then $J_1\cap V$ and $J_2\cap V$ are disjoint compact totally disconnected subsets of V.  It follows that there is a disk $D$ in $V$ such that $J_1\cap V\subset int(D)$ and $(J_2\cap V)\cap D=\emptyset$.  To see this, let $D^\prime$ be a disk in $V$, let $A_1$ be a subset of $int(D^\prime)$ that is homeomorphic to $J_1\cap V$, and let $A_2$ be a subset of $V-D^\prime$ that is homeomorphic to $J_2\cap V$.  Then according to Theorem 13.7 on pages 93-95 of   \cite{moise2013geometric}, there is a homeomorphism $\phi:V\rightarrow V$ such that $\phi(A_i)=J_i\cap V$ for $i=1,2$.  Simply let $D=\phi(D^\prime)$.  

We can assume $D$ is a piecewise linear disk in $V$.  Let $U$ be a regular neighborhood of $\partial D$ in $S^3$ such that $U\cap V$ is a regular neighborhood of $\partial D$ in $V$ and $U\cap J=\emptyset$, and let $T=cl(S^3-U)$.  
Then $U$ and, hence, $T$ are unknotted solid tori in $S^3$, and $J\subset int(T)$.  Let $\bar V=V\cup\{\infty\}$, and let $E_1$ and $E_2$ be the components of $cl(\bar V-U)$ such that $E_1\subset int(D)$ and $\infty\in E_2$.  Then $E_1$ and $E_2$ are disjoint meridional disks of $T$ such that $T\cap\bar V=E_1\cup E_2$, $J_1\cap\bar V\subset int(E_1)$ and $J_2\cap\bar V\subset int(E_2)$.  Thus, $J_1\cap E_2=\emptyset=J_2\cap E_1$. 

Clearly, there is a simple closed curve $K\subset int(T)$ such that the inclusion $K\hookrightarrow T$ is a homotopy equivalence, $p$ and $q\in K$, $K_1$ and $K_2$ are arcs such that $K_1\cup K_2=K$, $K_1\cap K_2=\partial K_1=\partial K_2=\{p,q\}$, and $K_1\cap E_2=\emptyset=K_2\cap E_1$.  At this point we will stretch conventional terminology slightly by saying that for subsets $Z\subset Y$ and $Z\subset Y^\prime$ of a space $X$, the inclusions $Y\hookrightarrow X$ and $Y^\prime\hookrightarrow X$ are \emph{homotopic in} $X$ \emph{rel} $Z$ if there is a homotopy $\xi:Y\times I\rightarrow X$ such that $\xi_0=id_Y$, $\xi_1:Y\rightarrow Y^\prime$ is a homeomorphism and $\xi_t\vert Z=id_Z$ for every $t\in I$.  Since $T-E_2$ is contractible, then the inclusions $J_1\hookrightarrow T-E_2$ and $K_1\hookrightarrow T-E_2$ are homotopic rel $\{p,q\}$.  Similarly, since $T-E_1$ is contractible, the inclusions $J_2\hookrightarrow T-E_1$ and $K_2\hookrightarrow T-E_1$ are homotopic rel $\{p,q\}$.  Therefore, the inclusions $J\hookrightarrow T$ and $K\hookrightarrow T$ are homotopic.  It follows that the inclusion $J\hookrightarrow T$ is a homotopy equivalence.

\textbf{Step 2:}  Consider a homeomorphism $\psi:B^2\times S^1\rightarrow T$, let $o\in int(B^2)$ and let $K=\psi(\{o\}\times S^1)$. Define the homeomorphism $\psi_o:S^1\rightarrow K$ by $\psi_o(z)=\psi(o,z)$.  Define the map $\tau:B^2\times S^1\rightarrow S^1$ by $\tau(y,z)=z$, and define the map $\pi:T\rightarrow K$ by $\pi=\psi_o\circ\tau\circ\psi^{-1}$.  Then we have a commutative diagram in which the vertical arrows are homeomorphisms:
\[
    \begin{tikzcd}
        B^2\times S^1   \arrow{rr}{\tau} 
        \arrow{d}[left]{\psi} && S^1 \arrow{d}{\psi_o} \\ T \arrow{rr}[below]{\pi} && K
    \end{tikzcd}
\]
We can now invoke Theorem 2 to obtain a level-preserving homeomorphism $\Theta:Cyl(\pi)\rightarrow Swl(\pi)$ that fixes the $0$- and $\infty$-levels.  Also observe that since $\tau:B^2\times S^1\rightarrow S^1$ is a homotopy equivalence, then so is $\pi:T\rightarrow K$.

\textbf{Step 3:}  Let $\lambda:I\rightarrow[0,\infty]$ be an order-preserving homeomorphism.  Clearly, there is an embedding $j:Cyl(\pi)\rightarrow S^3\times I$ with the following properties.\\
1)  $j$ maps $L(\pi,0)$ “identically” onto $T\times\{0\}$; i.e., $j(0x,\pi(x))=(x,0)$ for every $x\in T$.\\
2)  For every $t\in I$, $j(L(\pi,\lambda(t))\subset S^3\times\{t\}$, and $j(L(\pi,\lambda(t))$ is a copy of $T$ that is “squeezed” toward $K$.\\
3)  j maps $L(\pi,\infty)$ “identically” onto $K\times\{1\}$; i.e., for every $y\in K$, $j(\{\infty T\}\times\{y\})=\{(y,1)\}$.\\ 
Thus, $j\circ\Theta^{-1}:Swl(\pi)\rightarrow S^3\times I$ is an embedding with the following properties.\\
1)  $j\circ\Theta^{-1}$ maps $\mathcal{L}(\pi,0)$ “identically” onto $T\times\{0\}$; i.e., $j\circ\Theta^{-1}(0x,\pi(x))=(x,0)$ for every $x\in T$.\\
2)  For every $t\in I$, $j\circ\Theta^{-1}(\mathcal{L}(\pi, \lambda(t))\subset S^3\times\{t\}$.\\
3)  $j\circ\Theta^{-1}$ maps $\mathcal{L}(\pi,\infty)$ “identically” onto $K\times\{1\}$; i.e., for every $y\in K$, $j\circ\Theta^{-1}({\infty T}\times\{y\})=\{(y,1)\}$. 

\textbf{Step 4:}  Since $\pi\vert J\rightarrow K$ is the composition of the inclusion $J\hookrightarrow T$ and $\pi:T\rightarrow K$, both of which are homotopy equivalences, then $\pi\vert J\rightarrow K$ is a homotopy equivalence.  Hence, $\pi\vert J\rightarrow K$ is homotopic to an homeomorphism $\chi:J \rightarrow K$.  Therefore, Theorem 1 provides a $0$- and $\infty$-level preserving homeomorphism $\omega:Swl(\pi\vert J)\rightarrow Swl(\chi)$.  Thus, $j\circ\Theta^{-1}:Swl(\chi)\rightarrow S^3\times I$ is an embedding that maps $\mathcal{L}(\chi,0)$ into $S^3\times\{0\}$ and maps $\mathcal{L}(\chi,1)$ onto $K\times\{1\}\subset S^3\times\{1\}$.

\textbf{Step 5:}  Since $\chi:J\rightarrow K$ is a homeomorphism, then Theorem 2 provides a level-preserving homeomorphism from $Cyl(\chi)$ to $Swl(\chi)$.  Also, since $\chi:J\rightarrow K$ is a homeomorphism, then $Cyl(\chi)$ is an annulus with boundary $L(\chi,0)\cup L(\chi,1)$.  Therefore, $Swl(\chi)$ is an annulus with boundary $\mathcal{L}(\chi,0)\cup\mathcal{L}(\chi,1)$.  Since $\omega:Swl(\pi\vert J)\rightarrow Swl(\chi)$ is a $0$- and $\infty$-level preserving homeomorphism, then if follows that $Swl(\pi\vert J)$ is an annulus with boundary $\mathcal{L}(\pi\vert J,0)\cup \mathcal{L}(\pi\vert J,\infty)$.  Let $A=j\circ\Theta^{-1}(Swl(\pi\vert J))$.  Then A is an annulus in $S^3\times I$ with boundary $j\circ\Theta^{-1}(\mathcal{L}(\pi\vert J,0)\cup \mathcal{L}(\pi\vert J,\infty))=(J\times\{0\})\cup(K\times\{1\})$.

\textbf{Step 6:} Recall that for every $t\in I$, $j\circ\Theta^{-1}(\mathcal{L}(\pi\vert J,\lambda(t))\subset S^3\times\{t\}$.  Proposition 3 implies that there is a homeomorphism  $\Lambda:J\times[0,1)\rightarrow Swl(\pi\vert J)-\mathcal{L}(\pi\vert J,\infty)$ such that  $\Lambda(J\times\{t\})=\mathcal{L}(\pi\vert J,\lambda(t))$  for every $t\in[0,1)$.  Therefore, $j\circ\Theta^{-1}\circ\Lambda:J\times[0,1)\rightarrow A-(K\times\{1\})$ is a homeomorphism such that $j\circ\Theta^{-1}\circ\Lambda(J\times\{t\})\subset S^3\times\{t\}$ for every $t\in[0,1)$. $\blacksquare$

\section{Proofs of Theorems 1 and 2 and Proposition 3}

\textbf{Proof of Theorem 1.}  Let $X$ be a compact metric space and let $f, g:X\rightarrow S^1$ be homotopic maps.

\textbf{Step 1:}  There is a homeomorphism $\Phi:(\Sigma X)\times S^1\rightarrow(\Sigma X)\times S^1$ that carries $DblSwl(f)$ onto $DblSwl(g)$.

For each $x\in X$, observe that $\mathcal{DblF}(f,x)\cup\mathcal{DblF}(g,x)$ is a \emph{double helix} in the \emph{cylinder} $(\mathbb{R}x)\times S^1\subset(\Sigma X)\times S^1$.  We will construct a map $\sigma:X\rightarrow\mathbb{R}$ such that for each $x\in X$, a first-coordinate shift of $(\mathbb{R}x)\times S^1$ through a distance of $-\sigma(x)$ slides $\mathcal{DblF}(f,x)$ onto $\mathcal{DblF}(g,x)$.  (See Figure 3.)

\begin{figure}
    \centering
    \includegraphics[scale=.7]{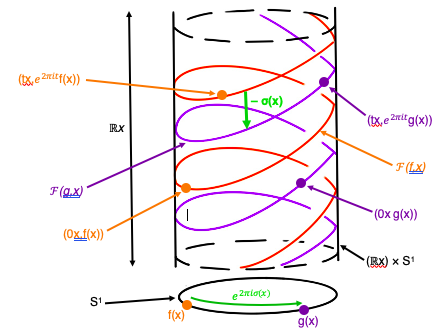}
    \caption{Double helix with vertical shift $-\sigma(x)$}
\end{figure}

Begin the construction of $\sigma$ by choosing a homotopy $h:X\times I\rightarrow S^1$ such that $h_0=f$ and $h_1=g$.  Using complex division in $S^1$, define $k:X\times I\rightarrow S^1$ by $k(x,t)=h(x,t)/h(x,0)$.  Then for every $x\in X$, $k(x,0) = 1$ and $k(x,1)f(x)=g(x)$.  Let $e:\mathbb{R}\rightarrow S^1$ be the exponential covering map $e(t)=e^{2\pi it}$.  The homotopy $k:X\times I\rightarrow S^1$ lifts to a homotopy $\widetilde{k}:X\times I\rightarrow \mathbb{R}$ such that $e\circ\widetilde{k}=k$ and $\widetilde{k}(x,0)=0$ for every $x\in X$.  Define $\sigma:X\rightarrow\mathbb{R}$ by $\sigma(x)=\widetilde{k}(x,1)$.  Therefore, for every $x\in X$,
\begin{center}
    $e^{2\pi i\sigma(x)}f(x)=e\circ\sigma(x)f(x)=e\circ\widetilde{k}(x,1)f(x)=k(x,1)f(x)=g(x)$.
\end{center}
Since X is compact, there is a $b>0$ such that $\sigma(X)\subset(-b,b)$.

Define $\Phi:(\Sigma X)\times S^1\rightarrow(\Sigma X)\times S^1$ by
\begin{equation*}
    \left\{
        \begin{array}{l}
            \Phi(tx,z)=((t-\sigma(x))x,z)\text{ for }(t,x)\in\mathbb{R}\times X \text{ and } z\in S^1, \text{ and}\\
            \Phi=id \text{ on } \{(-\infty)X,\infty X\}\times S^1.
        \end{array}
    \right.
\end{equation*}
$\Phi$ is the first-coordinate shift of $(\mathbb{R}x)\times S^1$ through a distance of $-\sigma(x)$ mentioned above.  It remains to show that $\Phi$ is a homeomorphism of $(\Sigma X)\times S^1$ which slides $\mathcal{DblF}(f,x)$ onto $\mathcal{DblF}(g,x)$ for each $x\in X$.

$\Phi$ is continuous at points of $\{(-\infty)X,\infty X\}\times S^1$ because for every $x\in X$ and every $z\in S^1$, $\Phi(([t,\infty]x)\times\{z\})\subset((t-b,\infty]x)\times\{z\}$ and $\Phi(([-\infty,t]x)\times\{z\})\subset([-\infty,t+b)x)\times\{z\}$.

$\Phi$ is a homeomorphism because its inverse $\overline{\Phi}$ can be defined explicitly by the equations
\begin{equation*}
    \left\{
        \begin{array}{l}
            \overline{\Phi}(tx,z)=((t+\sigma(x))x,z)\text{ for }(t,x)\in\mathbb{R}\times X \text{ and } z\in S^1, \text{ and}\\
            \overline{\Phi}=id \text{ on } \{(-\infty)X,\infty X\}\times S^1.
        \end{array}
    \right.
\end{equation*}
The verification that $\overline{\Phi}\circ\Phi=id=\Phi\circ\overline{\Phi}$ is straightforward.

For each $x\in X$, a typical point of $\mathcal{DblF}(f,x)$ has the form $(tx,e^{2\pi it}f(x))$, and
\begin{equation*}
    \begin{array}{l}
        \Phi((tx,e^{2\pi it}f(x))=((t-\sigma(x))x,e^{2\pi it}f(x))=\\
        ((t-\sigma(x))x,e^{2\pi i(t-\sigma(x))}e^{2\pi i\sigma(x)}f(x))=\\
        ((t-\sigma(x))x,e^{2\pi i(t-\sigma(x))}g(x))\in\mathcal{DblF}(g,x).
    \end{array}
\end{equation*}
Hence, $\Phi(\mathcal{DblF}(f,x))\subset\mathcal{DblF}(g,x)$.  A similar calculation shows 
\begin{equation*}
    \overline{\Phi}(\mathcal{DblF}(g,x))\subset\mathcal{DblF}(f,x) 
\end{equation*}
for every $x\in X$.  Thus,
\begin{equation*}
    \mathcal{DblF}(g,x)\subset\Phi(\mathcal{DblF}(f,x))
\end{equation*}
for every $x\in X$.  We conclude that
\begin{equation*}
    \Phi(\mathcal{DblF}(f,x))=\mathcal{DblF}(g,x) 
\end{equation*}
for every $x\in X$.

We have shown that $\Phi:(\Sigma X)\times S^1\rightarrow (\Sigma X)\times S^1$ is a homeomorphism with the property that $\Phi\vert DblSwl(f):DblSwl(f)\rightarrow DblSwl(g)$ is a fiber-preserving homeomorphism that fixes the $(-\infty)$- and $\infty$-levels.  Unfortunately, we can't conclude that $\Phi(Swl(f))$ and $Swl(g)$ are equal subsets of $DblSwl(g)$.  This completes Step 1.

\textbf{Step 2:}  There is a homeomorphism $\Psi:(\Sigma X)\times S^1\rightarrow(\Sigma X)\times S^1$ that maps $DblSwl(g)$ onto itself and \emph{twists} $\Phi(Swl(f)$ onto $Swl(g)$.

Observe that for every $x\in X$, $\Phi(\mathcal{F}(f,x))$ and $\mathcal{F}(g,x)$ are contained in the \emph{helix} $\mathcal{DblF}(g,x)$ which lies in the \emph{cylinder} $(\mathbb{R}x)\times S^1$.  For each $x\in X$, $\Psi$ will map this cylinder to itself with a \emph{screw motion} that preserves $\mathcal{DblF}(g,x)$ and \emph{twists} $\Phi(\mathcal{F}(f,x))$ onto $\mathcal{F}(g,x)$. 

Recall that $\sigma(X)\subset(-b,b)$.  Hence, there is clearly a map $\tau:\mathbb{R}\times X\rightarrow\mathbb{R}$ such that for each $x\in X$, $\tau\vert\mathbb{R}\times\{x\}:\mathbb{R}\times\{x\}\rightarrow\mathbb{R}$ is a homeomorphism such that $\tau(-\sigma(x))=0$ and $\tau(t,x)=t$ for $t\in(-\infty,-b]\cup[b,\infty)$.

Define $\Psi:(\Sigma X)\times S^1\rightarrow(\Sigma X)\times S^1$ by
\begin{equation*}
    \left\{
        \begin{array}{l}
            \Psi(tx,z)=((\tau(t,x)x,e^{2\pi i(\tau(t,x)-t)}z)\text{ for }(t,x)\in\mathbb{R}\times X \text{ and } z\in S^1, \text{ and}\\
            \Psi=id \text{ on } \{(-\infty)X,\infty X\}\times S^1.
        \end{array}
    \right.
\end{equation*}
$\Psi$ is continuous at points of $\{(-\infty)X,\infty X\}\times S^1$ because $\Psi=id$ on $\{tx:|t|\geq b\text{ and }x\in X\}\times S^1$.

$\Psi$ is a homeomorphism because its inverse $\overline{\Psi}$ can be defined explicitly as follows.  First observe that there is a map $\overline{\tau}:\mathbb{R}\times X\rightarrow \mathbb{R}$ so that for every $x\in X$, $t\mapsto \overline{\tau}(t,x):\mathbb{R}\rightarrow\mathbb{R}$ is the inverse of the homeomorphism $t\mapsto\tau(t,x):\mathbb{R}\rightarrow\mathbb{R}$.  Now define $\overline{\Psi}:(\Sigma X)\times S^1\rightarrow(\Sigma X)\times S^1$ by
\begin{equation*}
    \left\{
        \begin{array}{l}
            \overline{\Psi}(tx,z)=((\overline{\tau}(t,x)x,e^{2\pi i(\overline{\tau}(t,x)-t)}z)\text{ for }(t,x)\in\mathbb{R}\times X \text{ and } z\in S^1, \text{ and}\\
            \overline{\Psi}=id \text{ on } \{(-\infty)X,\infty X\}\times S^1.
        \end{array}
    \right.
\end{equation*}
The verification that $\overline{\Psi}\circ\Psi=id=\Psi\circ\overline{\Psi}$ is straightforward.

Let $x\in X$.  To prove that $\Psi(\mathcal{DblF}(g,x))=\mathcal{DblF}(g,x)$, note that a typical point of $\mathcal{DblF}(g,x)$ has the form $(tx,e^{2\pi it}g(x))$, and 
\begin{equation*}
    \begin{array}{l}
        \Psi((tx,e^{2\pi it}g(x))=(\tau(t,x)x,e^{2\pi i(\tau(t,x)-t)}e^{2\pi it}g(x))=\\
        (\tau(t,x)x,e^{2\pi i\tau(t,x)}g(x))\in\mathcal{DblF}(g,x).
    \end{array}
\end{equation*}
Therefore, $\Psi(\mathcal{DblF}(g,x))\subset\mathcal{DblF}(g,x)$.  A similar calculation shows 
\begin{equation*}
    \overline{\Psi}(\mathcal{DblF}(g,x))\subset\mathcal{DblF}(g,x). 
\end{equation*}
Hence,
\begin{equation*}
    \mathcal{DblF}(g,x)\subset\Psi(\mathcal{DblF}(g,x)).
\end{equation*}
We conclude that
\begin{equation*}
    \Psi(\mathcal{DblF}(g,x))=\mathcal{DblF}(g,x). 
\end{equation*}
We have shown that $\Psi:(\Sigma X)\times S^1\rightarrow (\Sigma X)\times S^1$ is a homeomorphism with the property that $\Psi\vert DblSwl(g):DblSwl(g)\rightarrow DblSwl(g)$ is a fiber-preserving homeomorphism that fixes the $(-\infty)$- and $\infty$-levels. 

Again let $x\in X$.  To prove that $\Psi(\Phi(\mathcal{F}(f,x)))=\mathcal{F}(g,x)$, note that a typical point of $\Phi(\mathcal{F}(f,x))$ has the form $(tx,e^{2\pi it}g(x)$ where $t\geq-\sigma(x)$, and a typical point of $\mathcal{F}(g,x)$ has the form $(tx,e^{2\pi it}g(x))$ where $t\geq0$.  Also note that $\tau$ maps $[-\sigma(x),\infty)$ onto $[0,\infty)$.  Since
\begin{equation*}
    \begin{array}{l}
        \Psi((tx,e^{2\pi it}g(x))=(\tau(t,x)x,e^{2\pi i(\tau(t,x)-t)}e^{2\pi it}g(x))=(\tau(t,x)x,e^{2\pi i\tau(t,x)}g(x)),
    \end{array}
\end{equation*}
then clearly $\Psi(\Phi(\mathcal{F}(f,x)))=\mathcal{F}(g,x)$.

It follows that $\Psi\circ\Phi(Swl(f))=Swl(g)$, completing Step 2.

We conclude that $\Psi\circ\Phi\vert Swl(f):Swl(f)\rightarrow Swl(g)$ is a fiber-preserving homeomorphism that fixes the $\infty$-level. $\blacksquare$

\textbf{Proof of Theorem 2.}  Let $Y$ be a compact metric space, $n$ a non-zero integer and $f:Y\times S^1\rightarrow S^1$ the map satisfying $f(y,z)=z^n$.  We will construct a homeomorphism $\Theta:C(Y\times S^1)\times S^1\rightarrow C(Y\times S^1)\times S^1$ that carries $Cyl(f)$ onto $Swl(f)$.

Define $\zeta:[0,\infty)\rightarrow S^1$ by $\zeta(t)=e^{-2\pi it/n}$, and define $\Theta:C(Y\times S^1)\times S^1\rightarrow C(Y\times S^1)\times S^1$ by
\begin{equation*}
    \left\{
        \begin{array}{l}
            \Theta(t(y,z),w)=(t(y,\zeta(t)z),w)\text{ for }t\in[0,\infty),(y,z)\in Y\times S^1\text{ and }w\in S^1,\text{ and}\\
            \Theta=id \text{ on }\{\infty(Y\times S^1)\}\times S^1.
        \end{array}
    \right.
\end{equation*}
$\Theta$ is clearly level-preserving and fixes the $0$- and $\infty$-levels.
 
 We argue that $\Theta$ is continuous at points of $\{\infty(Y\times S^1)\}\times S^1$.  For $a>0$, let $V_a=\bigcup_{t\in(a,\infty]}t(Y\times S^1)$; and for $w\in S^1$, let $\mathcal{M}_w$ be a basis for the topology on $S^1$ at $w$.  Then for each $w\in S^1$, $\{V_a\times M:a>0\text{ and }M\in\mathcal{M}_w\}$ is a basis for the topology on $C(Y\times S^1)\times S^1$ at $(\infty(Y\times S^1),w)$.  Since for each $a>0$ and each $M\in\mathcal{M}_w$, $\Theta$ maps $V_a\times M$ into itself, then $\Theta$ is continuous at $(\infty(Y\times S^1),w)$.
        
$\Theta$ is a homeomorphism because its inverse $\overline{\Theta}$ can be defined explicitly as follows.  First define the map $\overline{\zeta}:[0,\infty)\rightarrow S^1$ by $\overline{\zeta}(t)=e^{2\pi it/n}$.  Then define $\overline{\Theta}:C(Y\times S^1)\times S^1\rightarrow C(Y\times S^1)\times S^1$ by
\begin{equation*}
    \left\{
        \begin{array}{l}
            \overline{\Theta}(t(y,z),w)=(t(y,\overline{\zeta}(t)z),w)\text{ for }t\in[0,\infty),(y,z)\in Y\times S^1\text{ and }w\in S^1,\text{ and}\\
            \overline{\Theta}=id\text{ on }\{\infty(Y\times S^1)\}\times S^1
        \end{array}
    \right.
\end{equation*}
The verification that $\overline{\Theta}\circ\Theta=id=\Theta\circ\overline{\Theta}$ is straightforward.

To prove that $\Theta(Cyl(f))=Swl(f)$, note that a typical point of $Cyl(f)$ has the form $(t(x,z),f(x,z))$, a typical point of $Swl(f)$ has the form $(t(x,z),e^{2\pi it}f(x,z))$, and
\begin{equation*}
    \begin{array}{l}
        \Theta(t(x,z),f(x,z))=(t(x,\zeta(t)z),f(x,z))=\\
        (t(x,\zeta(t)z),z^n)=(t(x,\zeta(t)z),(\zeta(t))^{-n}(\zeta(t)z)^n)=\\
        (t(x,\zeta(t)z),e^{2\pi it}(\zeta(t)z)^n)=
        (t(x,\zeta(t)z),e^{2\pi it}f(x,\zeta(t)z))\in Swl(f).
    \end{array}
\end{equation*}
Therefore, $\Theta(Cyl(f))\subset Swl(f)$.  A similar calculation shows $\overline{\Theta}(Swl(f))\subset Cyl(f)$.  Hence, $Swl(f)\subset\Theta(Cyl(f))$.  We conclude that $\Theta(Cyl(f))=Swl(f)$.

It follows that $\Theta\vert Cyl(f):Cyl(f)\rightarrow Swl(f)$ is a level-preserving homeomorphism that fixes the $0$- and $\infty$-levels. $\blacksquare$

\textbf{Proof of Proposition 3.}  A compact metric space $X$, a map $f:X\rightarrow S^1$ and a homeomorphism $\lambda:[0,1)\rightarrow[0,\infty)$ are given.  Define $\Lambda:X\times[0,1)\rightarrow Swl(f)-\mathcal{L}(f,\infty)$ by $\Lambda(x,t)=(\lambda(t)x,e^{2\pi i\lambda(t)}f(x))$.  We show that $\Lambda$ is a homeomorphism by exhibiting its inverse.  Let $q:X\times[0,\infty]\rightarrow CX$ denote the quotient map $q(x,t)=tx$. Then $q\vert X\times[0,\infty):X\times[0,\infty)\rightarrow CX-\{\infty X\}$ is a homeomorphism; let $r:CX-\{\infty X\}\rightarrow X\times[0,\infty)$ denote its inverse.  Let $p:(CX-\{\infty X\})\times S^1\rightarrow CX-\{\infty X\}$ denote projection.  Define $\tau:(CX-\{\infty X\})\times S^1\rightarrow X\times[0,1)$ by $\tau=(id_X\times\lambda^{-1})\circ r\circ
p$.  It is easily verified that $\tau\circ\Lambda=id_{X\times[0,1)}$ and $\Lambda\circ(\tau\vert(Swl(f)-\mathcal{L}(f,\infty)))=id_{Swl(f)-\mathcal{L}(f,\infty)}$.  Hence, $\tau\vert(Swl(f)-\mathcal{L}(f,\infty))$ is the inverse of $\Lambda$. $\blacksquare$

\printbibliography

@article{ancel1996mapping,
  title={Mapping swirls and pseudo-spines of compact 4-manifolds},
  author={Ancel, Fredric D and Guilbault, Craig R},
  journal={Topology and its Applications},
  volume={71},
  number={3},
  pages={277--293},
  year={1996},
  publisher={Elsevier}
}

@article{bing1956simple,
  title={A simple closed curve that pierces no disk},
  author={Bing, R H},
  journal={J. Math. Pures Appl.(9)},
  volume={35},
  pages={337--343},
  year={1956}
}

@book{daverman1986decompositions,
  title={Decompositions of manifolds},
  author={Daverman, Robert J},
  year={1986},
  publisher={Academic Press}
}

@book{hillman2012algebraic,
  title={Algebraic invariants of links},
  author={Hillman, Jonathan Arthur},
  year={2012},
  publisher={World Scientific}
}

@incollection{matsumoto1979wild,
  title={Wild embeddings of piecewise linear manifolds in codimension two},
  author={Matsumoto, Yukio},
  booktitle={Geometric Topology},
  pages={393--428},
  year={1979},
  publisher={Elsevier}
}

@article{mazur1961note,
  title={A note on some contractible 4-manifolds},
  author={Mazur, Barry},
  journal={Annals of Mathematics},
  pages={221--228},
  year={1961},
  publisher={JSTOR}
}

@article{melikhov2020topological,
  title={Topological isotopy and Cochran's derived invariants},
  author={Melikhov, Sergey A},
  journal={arXiv preprint arXiv:2011.01409},
  year={2020}
}

@book{moise2013geometric,
  title={Geometric topology in dimensions 2 and 3},
  author={Moise, Edwin E},
  volume={47},
  year={2013},
  publisher={Springer Science \& Business Media}
}



\end{document}